\documentclass[amsthm]{elsart}

\usepackage{yjsco}
\usepackage{natbib}
\usepackage{amsthm}
\usepackage{enumerate}
\usepackage{amsmath}
\usepackage{amssymb}
\usepackage{amsfonts}
\usepackage{epsfig}
\usepackage{color,graphicx}

\newcommand{\qp}{\mathbb{Q}}
\newcommand{\pp}{\mathbb{P}}

\newcommand{\kp}{k}

\newcommand{\wt}{\widetilde}
\DeclareMathOperator{\coef}{coeff}

\DeclareMathOperator{\res}{res}

\begin{document}
\begin{frontmatter}

\title{Factoring bivariate polynomials using adjoints}

\author{Martin Weimann}
\address{Ricam, Austrian Academy of Sciences, Altenbergerstrasse 69,
A-4040 Linz, Austria}
\ead{weimann23@gmail.com}
\ead[url]{http://people.ricam.oeaw.ac.at/m.weimann/}

\begin{abstract}
One relates factorization of bivariate polynomials to singularities of projective plane curves. 
One proves that adjoint polynomials permit to find the recombinations of factors modulo $(x)$ induced by both absolute and rational factorizations, and so without using Hensel's lifting. One establishes in such a way the relations between the algorithms of Duval-Ragot based on locally constant functions and the algorithms of Ch\`eze-Lecerf based on lifting and recombinations. One shows in such a way that a fast computation of adjoint polynomials leads to a fast factorization.
The proof is based on cohomological sequences and residue theory.
\end{abstract}


\end{frontmatter}

\section{Introduction}

Factorization of multivariate polynomials is a central topic in Computer Algebra. One refers to \citep{C, CL, GG} and to the references therein for recent surveys of the topic. In this article, one studies the relations between singularities of projective plane curves and factorization of bivariate polynomials by using adjoint polynomials. In \citep{duval}, an algorithm is given for absolute bivariate factorization based on locally constant rational functions on the curve, using normalization and rational Newton-Puiseux expansions. The best actual complexities for rational \citep{L} and absolute \citep{CL} bivariate factorization have been obtained later on, based on a method of lifting and recombination of modular factors. One establishes here the bridge between these two approaches and one shows that the factorization can be computed fast from adjoint polynomials.  

\vskip4mm
\noindent
\textbf{Main result.} Let $F\in k[x,y]$ be a bivariate polynomial defined over a field $k$. 
We are interested in computing both the rational (over $k$) and absolute (over an algebraic closure $\bar{k}$) factorizations of $F$. In all of the sequel (except Section \ref{S7}), one assumes that $F$ satisfies the following hypothesis
\begin{equation*}
{\rm (H)} \qquad {\rm F(0,y) \,\,is \,\,separable\,\, of\,\, degree \,\,d=deg(F).}
\end{equation*}
In particular $F$ is square-free. Let $\mathcal{C}\subset \pp^2$ be the (reduced)	projective curve over $k$ defined by $F$.  One says that $D\subset \pp^2$ is an \textit{adjoint curve} of $\mathcal{C}$ if it passes throw all singular points $p$ of $\mathcal{C}$ (including infinitely near points) with multiplicity at least that of $\mathcal{C}$ minus one (see Section \ref{S4} for a more precise definition).
One says that $H\in k[x,y]$ is an \textit{adjoint polynomial} of $F$ of degree $n$ if it gives the dehomogenized equation of an adjoint curve of degree $n$. Adjoints may be computed by linear algebra from the resolution of singularities. One denotes by 
$$
A\subset k[y]
$$ 
the vector subspace spanned by the remainders modulo $(x)$ of adjoint polynomials of $F$ of degree $d-2$. Our main results assert that one can compute quickly both the rational and absolute factorizations of $F$ from the knowledge of a basis of $A$. 
\vspace{0.2cm} 
\noindent

One assumes that fast Fourier transform can be used for polynomial multiplication, so that two univariate polynomials over $k$ of degree $\le m$ can be multiplied in softly linear time $\wt{\mathcal{O}}(m)$. One denotes by $2\le\omega<3$ the matrix multiplication complexity exponent. 

\vspace{0.4cm} 
\noindent
\begin{thm}\label{t1} 
There exists a deterministic algorithm that, given $F$ satisfying (H) and given a basis of $A$, computes 
the rational factorization of $F$ with one factorization in $k[y]$ of degree $d$ plus
$$
\mathcal{O}(d^2(d-s)^{\omega-2})\subset \mathcal{O}(d^{\omega})
$$ 
arithmetic operations over $k$, with $s$ the number of irreducible rational factors of $F$.
\end{thm}

\vspace{0.4cm} 
\noindent
The vector space $A$ contains enough information to compute too the absolute factorization of $F$. By absolute factorization, one means here the computation of a family of pairs of polynomials
$$
\{(Q_1,q_1),\ldots,(Q_r,q_r)\}
$$

\vspace{0.2cm} 
\noindent
where $q_i\in k[t]$ is monic, $Q_i\in k[x,y,t]$ and where
$$
F(x,y)=\prod_{i=1}^r \prod_{q_i(\alpha)=0} Q_i(x,y,\alpha)
$$
is the irreducible decomposition of $F$ over $\bar{k}$. Note that such a representation is non unique. One obtains the following result :

\vspace{0.4cm} 
\noindent

\begin{thm}\label{absolute} Suppose that $k$ has characteristic $0$ or greater than $d(d-1)$. There exists a deterministic algorithm that, given $F$ satisfying hypothesis (H) and given a basis of $A$, computes the absolute factorization of $F$ within
$$
\wt{\mathcal{O}}(d^2(d-\bar{s})^{\omega-2}+\bar{s}d^{3})\subset \wt{\mathcal{O}}(d^4)
$$ 
arithmetic operations over $k$, with $\bar{s}$ the number of irreducible absolute factors.
\end{thm}
\vspace{0.4cm} 
\noindent
Note that in contrast to the rational case, no univariate factorization is required for absolute factorization. 
Following \citep{CL}, one obtains too a probabilistic approach in terms of computation trees. 

\vspace{0.4cm} 
\noindent
\begin{thm}\label{absoluteproba} 
Suppose that $k$ has characteristic $0$ or greater than $d(d-1)$. Given $F$ satisfying hypothesis (H) and given a basis of $A$, there exists a polynomial $S\in\bar{k} [t_1,\ldots,t_{d}]$ of degree at most $d(d-1)$ and a family of computation trees parametrized by 
$c\in \kp^d$ such that 

\vspace{0.1cm} 
\noindent

\hspace{0.5cm} $\bullet$ Any executable tree returns the absolute factorization of $F$;

\hspace{0.5cm} $\bullet$ A tree is executable whenever $S(c)\ne 0$. 

\vspace{0.1cm} 
\noindent
The maximal cost of the trees is bounded by
$$
\wt{\mathcal{O}}(d^2(d-\bar{s})^{\omega-2}+d^{\frac{\omega+3}{2}})\subset \wt{\mathcal{O}}(d^{\frac{\omega+3}{2}})
$$ 
arithmetic operations over $k$. If the cardinality of $\kp$ is infinite, the algorithm returns the correct answer with probability one.
\end{thm}

\vspace{0.4cm} 

Up to our knowledge, all our complexities are smaller or equal to the actual best complexities for factorization of dense bivariate polynomials that are known to be $\mathcal{O}(d^{\omega+1})$ for rational deterministic \citep{L}, $\wt{\mathcal{O}}(d^4)$ for absolute deterministic and $\wt{\mathcal{O}}(d^{3})$ for absolute probabilistic \citep{CL}. Thus, our results lead immediately to the following question :

\vskip3mm 
\noindent
\textbf{Question\label{q1}:} Can we compute a basis of $A$ fast enough for our method being useful for rational or absolute factorizations ?
\vskip3mm  
\noindent
One obtains the following encouraging result, as a consequence of the Riemann-Roch theorem for reducible curves.

\vspace{0.4cm} 

\begin{thm}\label{adjmodx} Given a basis of the vector space of all adjoint polynomials of degree $d-2$, one computes a basis of $A$ within 
$$
\mathcal{O}(d(g+d-\bar{s})(d-\bar{s})^{\omega-2})\subset \mathcal{O}(d^{\omega+1})
$$
arithmetic operations over $k$, where $g$ is the geometric genus of $\mathcal{C}$.
\end{thm}

\vskip4mm  
\noindent

Note that there exist efficient algorithms for computing adjoint polynomials, by using Newton-Puiseux expansions \citep{sw} or integral basis \citep{Mnuk, VH}, but whose complexities have not been analyzed yet. Unfortunately, it is {\it a priori} hopeless that Theorem \ref{adjmodx} answers positively to our question : one needs the all resolution of singularities to compute adjoints and 
one expects factorization to be one step in desingularization (up to our knowledge, if $k=\mathbb{F}_p$ with $p>d$, the best complexity for computing all singular Puiseux expansions of $F$ is $\wt{\mathcal{O}}(d^5)$ operations \citep{PR}). However, one will show (Section \ref{S8}) that it's enough to separate all local branches $\mathcal{C}$ in order to compute $A$. In particular, one needs not to desingularize irreducible branches. Moreover, one will see (Subsection \ref{S8}) that in some cases, our method adapts too to the case $F(0,y)$ non separable, in which case one can use the combinatorial information given by the resolution of singularities along the line $x=0$ to speed-up the algorithm. In that spirit, the author recently developed a factorization algorithm based on the toric resolution of the singularities at infinity \citep{W2}, running in polynomial time in the volume of the Newton polytope, improving \citep{L} for sparse enough polynomials. 

To summarize, the method developed here may be useful for some special type of polynomials (smooth components, high singularities along a line, etc.) or more generally if one is given some extra input data concerning the singularities. Anyway, our results clarify the relations between normalization and factorization, giving a good point of view for comparing the classes of complexity of both operations and of various related algorithms (Newton-Puiseux, Hensel lifting, integral closure, etc.). It has to be noticed too that the strength of our approach depends strongly on further improvements in the algorithmic theory of singularities, especially on the Newton-Puiseux algorithm. 

\vspace{0.3cm} 
The proofs of our main results rely on the structural sheaf sequence of a divisor on the normalization of $\mathcal{C}$, combined with the Serre duality and with the residue theorem. Roughly speaking, our algorithms combine ideas developed in \citep{duval} and \citep{ragot} (computing locally constant rational functions) with ideas in \citep{L} and \citep{CL} (lifting and recombination of modular factors). Namely, one shows that one can recombine the factors modulo $(x)$ from $A$ without using Hensel lifting. In fact, one proves that \textit{absolute recombinations and adjoints modulo $(x)$ determine each other by solving a $d\times d$ linear system over $k$} (see Corollary \ref{c2} for a precise statement).                                                                                                                       
\vskip3mm  
\noindent
\textbf{Organization.} One introduces the recombination problem and its relation to locally constant functions in Section \ref{S2}. In Section \ref{S3}, one proves our key result that gives conditions for lifting locally constant functions using residue theory and cohomology. In Section \ref{S4}, one establishes the relation with adjoint polynomials and proves Theorem \ref{adjmodx}. One solves recombinations in Section \ref{S5} from which follow the proofs of Theorem \ref{t1}, \ref{absolute} and \ref{absoluteproba} in Section \ref{S6}. In Section \ref{S7}, one discusses the case $F(0,y)$ non separable and one illustrates our approach on a simple example. In Section \ref{S8}, one shows that the computation of $A$ does not require the all resolution of singularities. Finally, one concludes in the last Section \ref{S9}.

\section{Recombinations and locally constant functions.}\label{S2}

Our algorithms are related to \citep{L}  and \citep{CL}, both methods being based on the recombination problem of the modulo $(x)$ factors.  One first explains this problem and then one relates it to the sheaves of locally constant functions on the normalizing curve. One keeps the same notations and hypothesis as in the introduction. %

\subsection{Recombinations problems}\label{ss2.1} Let us consider the respective factorizations 
$$
\begin{cases}
F(x,y) = F_1(x,y)\cdots F_s(x,y)\\
F(0,y) = f_1(y)\cdots f_n(y)
\end{cases}
$$
of $F$ and $F$ modulo $(x)$ over $k$ (recall that $F(0,y)$ is assumed to be separable). Solving {\it rational recombinations} consists in computing the vectors $$\nu^{(j)}=(\nu^{(j)}_{1},\ldots,\nu^{(j)}_n)\in \{0,1\}^n$$ induced by the relations
$$
F_j(0,y)=\prod_{i=1}^n f_i(y)^{\nu^{(j)}_i}, \,\,\, j=1,\ldots,s.
$$
In the same way, let
$$
\begin{cases}
F(x,y) = \bar{F}_1(x,y)\cdots \bar{F}_{\bar{s}}(x,y)\\
F(0,y) = \bar{f}_1(y)\cdots \bar{f}_d(y)
\end{cases}
$$
be the respective factorizations of $F$ and $F$ modulo $(x)$ over $\bar{k}$. Solving {\it absolute recombinations} consists in computing the vectors $$\bar{\nu}^{(j)}=(\bar{\nu}^{(j)}_1,\ldots,\bar{\nu}^{(j)}_d)\in \{0,1\}^d$$ induced by the relations
$$
\bar{F}_j(0,y)=\prod_{i=1}^d \bar{f}_i(y)^{\bar{\nu}^{(j)}_i}, \,\,\, j=1,\ldots,\bar{s}.
$$
The following picture illustrates the absolute recombinations when $\mathcal{C}$ is union of a cubic and a conic. 
\vskip3mm  
\noindent
\begin{equation*}
\begin{array}{c}
\resizebox{6.cm}{!}{\includegraphics{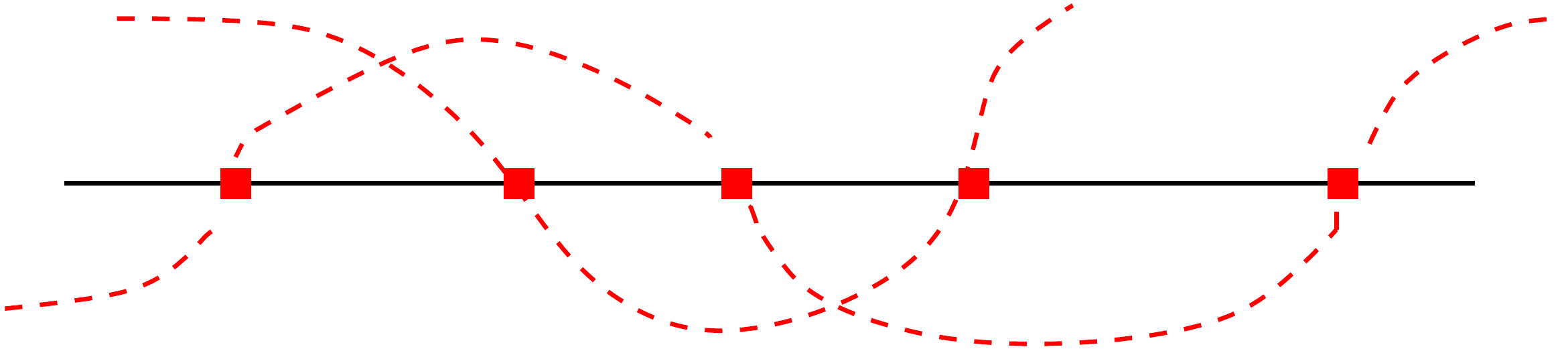}} 
\end{array}
\end{equation*}

$$
\Longrightarrow \qquad \bar{\nu}^{(1)}=(1,0,1,0,1),\quad \bar{\nu}^{(2)}=(0,1,0,1,0)\qquad\qquad
$$
\vskip3mm  
\noindent

In this article we mainly pay attention to the recombination problems, the irreducible factorization of $F$ then following with a fast multi-factor Hensel lifting (combined with a partial fraction decomposition algorithm in the absolute case). The main idea is to interpret the recombination problem as a cohomological problem of lifting sections. 

\vskip3mm  
\noindent

\subsection{Solving recombinations {\it {\bf via}} lifting sections}\label{ss2.2} All schemes and properties (connectivity, irreducibility) are considered over the base field $k$. One may think a point over a scheme $X$ over a $k$ as a collection of points in the extension $\bar{X}=X\otimes_k \bar{k}$ that are conjugated under the Galois group of $\bar{k}/k$. 

Let $\mathcal{C}$ and  $\mathcal{L}$ be the respective Zariski closures of the affine curves $F=0$ and $x=0$ to the projective plane $\pp^2$. Let
$$
\pi : X \rightarrow \pp^2
$$
be the standard embedded resolution of $\mathcal{C}$. One denotes by $C$ and $L$ the respective strict transforms of $\mathcal{C}$ and $\mathcal{L}$ by $\pi$. The inclusion of the zero-dimensional subscheme 
$$
Z:=C\cap L
$$ into $C$ induces a restriction morphism
$$
\alpha \,: \,H^0(\mathcal{O}_{C}) \hookrightarrow  H^0(\mathcal{O}_{Z}) 
$$
between the respective $k$-vector spaces of regular functions on $C$ and $Z$. Note that both vector spaces may be identified with the sets of locally constant functions on $C$ and $Z$. The map $\alpha$ is injective since $Z$ has at least one point on each component of $C$.  The following two subsections are dedicated to show that the computation of the cokernel of $\alpha$ permits to solve both the rational and the absolute recombination problems.

\vskip3mm  
\noindent

\subsubsection{The rational case.} 

\vskip2mm  
\noindent
Since $F(0,y)$ has degree $d$, $Z$ is an affine zero-dimensional subscheme whose ring of regular functions may be identified with the finite $k$-algebra
\begin{equation}\label{Z1}
H^0(\mathcal{O}_{Z}) = \frac{k[x,y]}{(x,F)} = \frac{k[y]}{(F(0,y))}.
\end{equation}
Since $F(0,y)$ is separable, its rational factorization induces an isomorphism
\begin{equation}\label{Z2}
H^0(\mathcal{O}_{Z}) \simeq \frac{k[y]}{(f_1)}\oplus\cdots\oplus \frac{k[y]}{(f_n)}.
\end{equation}
Thus $Z$ has $n$ connected components (closed points) $p_1,\ldots,p_n$ corresponding to the maximal ideals of the ring $H^0(\mathcal{O}_{Z})$ generated by the $f_i$'s. The natural inclusions 
$$
k\hookrightarrow \frac{k[y]}{(f_i)},\quad i=1,\ldots,n
$$ 
combined with (\ref{Z1}) and (\ref{Z2}) induce the inclusion  
$$
k^n\subset H^0(\mathcal{O}_{Z}),
$$ 
$k^n$ being identified with the subspace of locally constant functions on $Z$ that take value in $k$, that is $(\nu_1,\ldots,\nu_n)\in k^n$ sends $p_i$ to $\nu_i$ (in general, a function on $Z$ takes values in the various residue fields $k[y]/(f_i)$). The map $\alpha$ introduced before is related to recombinations by the following lemma:

\vskip3mm  
\noindent
\begin{lem}\label{lemW}
The vector subspace $W\subset k^n$ defined by 
$$
W:=k^n\cap Im(\alpha)
$$
admits $(\nu^{(1)},\ldots,\nu^{(s)})$ as reduced echelon basis (up to reordering). 
\end{lem}
\vskip2mm  
\noindent
\begin{proof}
By definition, $\nu\in W$  if and only if it's the restriction to $Z$ of a locally constant $k$-valued function on $C$. Since $C$ is smooth, it has $s$ connected components $C_1\ldots, C_s$ corresponding to the prime rational factors of $F$. 
Thus $\nu\in W$ if and only if $\nu$ is $k$-valued and constant along $C_j\cap L$ for $j=1,\ldots,s$. One deduces that $dim_{k} \, W = s$ and that $\nu^{(j)}\in W$ for $j=1,\ldots,s$.  
Since the $\nu^{(j)}$'s have $\{0,1\}$-coordinates and are pairwise orthogonal vectors in $k^n$, they form up to reordering the reduced echelon basis of $W$.
\end{proof}

By Lemma \ref{lemW}, the recombination problem over $k$ is reduced to compute first the rational factorization of $F(0,y)$ (inducing the inclusion $k^n\subset H^0(\mathcal{O}_{Z})$), and then the cokernel of $\alpha$. 

\vskip3mm  
\noindent
\subsubsection{The absolute case.} 
\vskip2mm  
\noindent

The relations between locally constant functions and absolute factorization is explored in \citep{duval} where the author determines one absolute factor a time from a basis of the regular functions on $C\times_{k}\bar{k}$. One rather relates here regular functions on $C$ to the recombination algorithm in \citep{CL} and one computes all irreducible factors simultaneously by using multi-factor Hensel lifting. One first proves :

\vskip2mm  
\noindent
\begin{lem}\label{lemDim}
One has equalities $dim_{k} H^0(\mathcal{O}_{Z})=d$ and $dim_{k} H^0(\mathcal{O}_{C})=\bar{s}$.
\end{lem}
\vskip2mm  
\noindent

\begin{proof}
First equality is clear from (\ref{Z1}). Since $H^0(\mathcal{O}_{C})$ is a finite dimensional $k$-vector space, one has 
$$
dim_{k} \, H^0(\mathcal{O}_{C})=dim_{\bar{k}} \,H^0(\mathcal{O}_{C})\otimes_{k}\bar{k}
$$
Let $\bar{C}:=C\times_{k}\bar{k}$ be the geometrical scheme associated to $C$ by extending the base field $k$ to its algebraic closure $\bar{k}$. One has \citep[prop. $1.24$ p.85]{Liu}
$$
H^0(\mathcal{O}_{C})\otimes_{k}\bar{k}=H^0(\mathcal{O}_{\bar{C}}).
$$
Since $\bar{C}$ is smooth, it's the \textit{disjoint} union of $\bar{s}$ irreducible components $\bar{C}_1\ldots, \bar{C}_{\bar{s}}$ corresponding in an obvious way to the prime absolute factors of $F$. It follows that one has an isomorphism of $\bar{k}$-vector spaces
$$
H^0(\mathcal{O}_{\bar{C}})\simeq \oplus_{j=1}^{\bar{s}} H^0(\mathcal{O}_{\bar{C}_j}).
$$
Since $H^0(\mathcal{O}_{\bar{C}_j})=\bar{k}$, one has $dim_{\bar{k}}H^0(\mathcal{O}_{\bar{C}})=\bar{s}$ so that $dim_{k} \, H^0(\mathcal{O}_{C})=\bar{s}$.
\end{proof}
\vskip3mm  
\noindent

Let $\phi_1,\ldots,\phi_d$ be the roots of $F(0,y)$ in $\bar{k}$. The identification (\ref{Z1}) gives rise to the multi-evaluation isomorphism

\begin{equation}\label{eval}
\begin{array}{ccccc}
& ev: & H^0(\mathcal{O}_{Z}) \otimes_{k}  \bar{k} &\stackrel{\simeq}{\longrightarrow} & \bar{k}^d \\
& & \nu & \longmapsto & \bar{\nu}:=(\nu(\phi_1),\ldots,\nu(\phi_d)).
\end{array}
\end{equation}
\vskip3mm  
\noindent
The next lemma shows that solving absolute recombinations reduces to compute $Im(\alpha)$ and to apply the evaluation map $ev$. One endows $\bar{k}^d$ with its canonical basis.

\vskip3mm  
\noindent
\begin{lem}\label{Wbar}
The vector subspace $\bar{W}\subset \bar{k}^d$ defined by
$$
\bar{W}:=ev(Im(\alpha)\otimes_{k} \bar{k})
$$
admits $(\bar{\nu}^{(1)},\ldots,\bar{\nu}^{(\bar{s})})$ as reduced echelon basis (up to reordering). 
\end{lem}
\vskip3mm  
\noindent

\begin{proof}
Let $\bar{Z}:=Z\times_{k}\bar{k}$. The map $ev$ induces an identification 
$$
\bar{k}^d=H^0(\mathcal{O}_{\bar{Z}}),
$$
where $\bar{\nu}=(\bar{\nu}_1,\ldots,\bar{\nu}_d)\in \bar{k}^d$ is identified with the locally constant function that sends each closed point $\bar{p_i}\in\bar{Z}$ to $\bar{\nu}_i$. Since $\bar{Z}$ contains at least one point of each connected component of $\bar{C}$, the restriction map
$$
\bar{\alpha}\,: \,H^0(\mathcal{O}_{\bar{C}}) \hookrightarrow  H^0(\mathcal{O}_{\bar{Z}})
$$ 
is injective. By definition, 
$
\bar{W}=Im(\bar{\alpha})
$
so that $dim_{\bar{k}} \, \bar{W} = dim_{\bar{k}} H^0(\mathcal{O}_{\bar{C}})=\bar{s}$ by the proof of Lemma \ref{lemDim}. Each vector $\bar{\nu}^{(j)}$ being constant on $\bar{C}_1\cap \bar{L},\ldots, \bar{C}_{\bar{s}}\cap \bar{L}$, it extends to a function on $\bar{C}$. So $\bar{\nu}^{(j)}\in \bar{W}$ for $j=1,\ldots,\bar{s}$. Since the $\bar{\nu}^{(j)}$'s have $\{0,1\}$-coordinates and are pairwise orthogonal in $\bar{k}^d$, they form up to reordering the reduced echelon basis of $\bar{W}$. 
\end{proof}

\vskip3mm  
\noindent

\section{Lifting sections using residues}\label{S3} 

The previous section shows that recombinations may be reduced to compute the cokernel of the restriction morphism
$$
\alpha \,: \,H^0(\mathcal{O}_{C}) \hookrightarrow  H^0(\mathcal{O}_{Z}).
$$
To this aim, one introduces residues. There is an extensive litterature \citep{couvreur, serre, lipman, tate, vakil} concerning residues on curves and surfaces. 

Let $\omega_C$ be the sheaf of regular differential $1$-forms over $C$ (the dualizing sheaf) and let $\omega_C(Z)$ be the sheaf of meromorphic $1$-forms with polar divisor bounded by $Z$. Let $p\in C$ with residue field $k_p$ and let $\psi\in\omega_{C,p}(Z)$ be a germ of meromorphic form at $p$. For any uniformizer $t$ of $C$ at $p$, there exists a unique formal series $h\in k_p[[t]]$ such that 
$$
\psi=\frac{h(t)dt}{t}.
$$
One defines \textit{the residue of $\psi$ at $p$} as 
$$
res_{p} \,\psi := Tr_{p}\, [h(0)],
$$
where $Tr_{p} : k_p \to k$ is the trace map. This definition does not depend on the choice of the uniformizer (see for instance \citep{serre}). The map $res_p$ is $k$-linear and vanishes on  regular forms. In particular, if  $\nu\in \mathcal{O}_{Z,p}$ has a local lifting $\wt{\nu}$ to $\mathcal{O}_{C,p}$, one checks that the definition
$$
\res_{p} (\nu\,\psi):=\res_{p} (\wt{\nu}\,\psi)
$$ 
does not depend on the choice of the lifting. One obtains the following key result.

\vskip4mm
\noindent

\begin{prop}\label{t2}
There is an exact sequence of $k$-vector spaces
\vskip2mm
\noindent
$$
0 \longrightarrow  H^0(\mathcal{O}_{C}) \stackrel{\alpha}\longrightarrow  H^0(\mathcal{O}_{Z}) \stackrel{R}{\longrightarrow} H^0(\omega_C(Z))^{\vee} \stackrel{\beta} 
\longrightarrow H^0(\omega_C)^{\vee} \longrightarrow 0
$$ 
\vskip2mm
\noindent
where $^{\vee}$ stands for the dual and where $R$ associates to $\nu$ the linear form 
\begin{equation*}
R_{\nu} :  \,\,\psi  \, \longmapsto  \, \sum_{i=1}^n \res_{p_i} (\nu\psi).
\end{equation*}
In particular, $dim\, H^0(\omega_C(Z))=g+d-\bar{s}$ where $g$ is the geometric genus of $\mathcal{C}$ (sum of the genus of the irreducible components).
\end{prop}
 
\vskip4mm
\noindent

\begin{proof}
Let $\omega_Z$ be the dualizing sheaf of $Z$. One has $\omega_Z\simeq Hom(\mathcal{O}_Z, k)$, and the local duality theorem gives a short exact sequence (the adjunction formula) 
$$
0 \longrightarrow \omega_C \longrightarrow \omega_C(Z)\stackrel{Res}\longrightarrow 
\omega_Z \longrightarrow 0,
$$
where the residue map $Res$ is locally defined on an open set $U\subset C$ as
\[\begin{aligned}
Res_U(\psi): \mathcal{O}_Z(U) & \longrightarrow  k \\
\nu & \longmapsto  \sum_{p \in U} res_{p}(\nu\psi).
\end{aligned}\]
The associated long exact cohomology sequence is

\begin{equation}\label{koszul}
0 \rightarrow H^0(\omega_C)\rightarrow H^0(\omega_C(Z)) \stackrel{Res}\rightarrow 
H^0(\omega_Z) \rightarrow H^1(\omega_C)\rightarrow H^1(\omega_C(Z)).
\end{equation}
\vskip4mm
\noindent
By the duality of Serre, one has isomorphisms 
$$H^1(\omega_C)\simeq H^0(\mathcal{O}_{C})^{\vee}\quad {\rm and} \quad H^1(\omega_C(Z))\simeq H^0(\mathcal{O}_{C}(-Z))^{\vee}=0,
$$
the last vanishing property because $Z$ as at least one point on each connected component of $C$. 
The dual sequence of (\ref{koszul}) becomes

$$
0 \rightarrow  H^0(\mathcal{O}_{C}) \stackrel{\alpha}\rightarrow  H^0(\mathcal{O}_Z) \stackrel{R}{\rightarrow} H^0(\omega_C(Z))^{\vee} \stackrel{\beta} 
\rightarrow H^0(\omega_C)^{\vee} \rightarrow 0
$$ 
\vskip4mm
\noindent
where $R$ is dual to $Res$, that is 
$$
R: \nu \longmapsto \Big(\psi\mapsto \sum_{p_i\in Z} res_{p_i}(\nu\psi)\Big).
$$
This shows the exact sequence of Proposition \ref{t2}.  This sequence induces equality
$$
h^0(\omega_C(Z))= h^0(\omega_C)+h^0(\mathcal{O}_Z)-h^0(\mathcal{O}_C)=g+d-\bar{s},
$$
last equality using Lemma \ref{lemDim} and using that $h^0(\omega_C)$ coincides with the geometric genus of $\mathcal{C}$.
This ends the proof.
\end{proof}

\vskip3mm
\noindent

\begin{rem}
The inclusion $Im(\alpha)\subset ker(R)$ follows from the \textit{residue theorem} that asserts that 
$$
\sum_{p\in C_j} res_p \psi=0
$$
for all connected component $C_j$ of $C$ and all rational $1$-form  $\psi$ on $C$.
\end{rem}

\vskip3mm
\noindent

\begin{rem}
In the case $C$ irreducible over $\bar{k}$, the equality $$dim\, H^0(\omega_C(Z))=g+d-1$$ given by Proposition \ref{t2} follows from the \textit{theorem of Riemann-Roch} for curves. 
\end{rem}

\vskip3mm
\noindent

\section{Relations with adjoint polynomials}\label{S4} 

One relates now holomorphic forms with adjoint polynomials. 
One denotes by $S$ the set of singular points of $\mathcal{C}$, including all infinitely near points. For each $p\in S$, there is a decomposition of $\pi$
$$
X\stackrel{\pi_1}\longrightarrow \wt{X}_p\stackrel{\pi_p}\longrightarrow X_p \stackrel{\pi_2} \longrightarrow\pp^2
$$
such that $p$ is a closed point of the intermediary surface $X_p$ and $\pi_p$ is the blow-up at $p$. Let $E_p$ be the exceptional divisor of $\pi_p$ and $\hat{E_p}$ its total transform under $\pi_1$.
One denotes by $m_p$ the multiplicity at $p$ of the strict transform of $\mathcal{C}$ under the map $\pi_2$. 

\vskip3mm
\noindent
\begin{defn}\label{defadj}
The \textit{adjoint divisor} of $F$ is the exceptional effective divisor 
$$
E:=\sum_{p\in S} (m_p-1)\hat{E_p}.
$$
An \textit{adjoint curve} of $\mathcal{C}$ is an effective divisor $D\subset \pp^2$ that satisfies
$$
\pi^*(D)\ge E.
$$
An \textit{adjoint polynomial} of $F$ of degree $\le m$ is a polynomial giving the dehomogeneised affine equation of an adjoint curve of degree $m$. 
\end{defn}

\vskip3mm
\noindent

In other words, adjoints of $F$ are those polynomials vanishing at the singular points of $\mathcal{C}$ with high enough multiplicities. Adjoints carry out precious informations about the geometry of $\mathcal{C}$. In particular, it is well known that they are deeply related to the sheaf $\omega_C$ of regular forms on the normalized curve. Let us denote by 
$$
Adj(m)\subset k[x,y]
$$ 
the $k$-vector subspace generated by adjoint polynomials of $F$ of degree $\le m$. One has the following proposition:

\vskip4mm
\noindent

\begin{prop}\label{p1}
For all integers $m\le 2$, one has an isomorphism
\[\begin{aligned}
Adj(d-3+m)  & \stackrel{\simeq}\longrightarrow  H^0(\omega_C(mZ))\\
H & \longmapsto  \pi^*\Big(\frac{Hdx}{x^m \partial_y F}\Big)_{|C}.
\end{aligned}\]
\end{prop}
\vskip2mm
\noindent
\begin{proof}
In order to relate adjoints with differential forms, one introduces the \textit{conductor} 
\vskip2mm
\noindent
$$
\mathcal{A}_{\mathcal{C}}:=Hom_{\mathcal{O}_{\mathcal{C}}}(\pi_*\mathcal{O}_C, \mathcal{O}_{\mathcal{C}})
$$ 
\vskip2mm
\noindent
of the normalization of $\mathcal{C}$. It is an ideal sheaf of $\mathcal{O}_{\mathcal{C}}$, related to the dualizing sheaf $\omega_{\mathcal{C}}$ of $\mathcal{C}$ by the formula
$$
\pi_*\omega_C=\omega_{\mathcal{C}}\otimes_{\mathcal{O}_{\mathcal{C}}} \mathcal{A}_{\mathcal{C}}
$$
\vskip2mm
\noindent
\citep[see for instance][p.25]{szpiro}. Since the morphism $\pi:C\to \mathcal{C}$ is affine one has 
\vskip3mm
\noindent
\begin{equation}\label{a1}
H^0(C,\omega_C(mZ))=H^0(\mathcal{C},\pi_*(\omega_C(mZ)))=H^0(\mathcal{C},\omega_{\mathcal{C}}(m\mathcal{L})\otimes \mathcal{A}_{\mathcal{C}}),
\end{equation}
\vskip3mm
\noindent
last equality following from the projection formula (recall that $Z=C\cap L$ and $\pi^* \mathcal{L}=L$). 
Let $\mathcal{A}$ be the inverse ideal sheaf of $\mathcal{A}_{\mathcal{C}}$ under the restriction $\mathcal{O}_{\pp^2}\to \mathcal{O}_{\mathcal{C}}$. One has the short exact sequence 
\vskip2mm
\noindent
\begin{equation}\label{a2}
0\longrightarrow \mathcal{O}_{\pp^2}(-\mathcal{C}) \longrightarrow \mathcal{A} \longrightarrow \mathcal{A}_{\mathcal{C}}\longrightarrow 0.
\end{equation}
\vskip4mm
\noindent
Tensoring (\ref{a2}) with the invertible sheaf $\Omega_{\pp^2}^2(\mathcal{C}+m\mathcal{L})$, and using the adjunction formula, one obtains the exact sequence
\vskip2mm
\noindent
\begin{equation}\label{a3}
0\longrightarrow \Omega_{\pp^2}^2(m\mathcal{L}) \longrightarrow \Omega_{\pp^2}^2(\mathcal{C}+m\mathcal{L})\otimes\mathcal{A} \stackrel{RP}\longrightarrow \omega_{\mathcal{C}}(m\mathcal{L})\otimes \mathcal{A}_{\mathcal{C}}\longrightarrow 0.
\end{equation}
\vskip2mm
\noindent
Here, $RP$ is the Poincar\'e residue map, defined outside the singular locus of $\mathcal{C}$ as
\vskip2mm
\noindent
\begin{equation}\label{a4}
RP\Big(\frac{Hdx\land dy}{F x^m}\Big)=\Big(\frac{Hdx}{x^m \partial_y F}\Big)_{|\mathcal{C}}.
\end{equation}
\vskip2mm
\noindent
For $m\le 2$ and $i=0,1$, one has \citep[Theorem 5.1 p.225]{hart}
\vskip2mm
\noindent
$$
H^i(\pp^2,\Omega^2_{\pp^2}(m\mathcal{L}))= H^i(\pp^2,\mathcal{O}_{\pp^2}(m-3))=0 
$$ 
\vskip2mm
\noindent
which, combined with the long exact cohomological sequence of (\ref{a3}), gives an isomorphism
\vskip2mm
\noindent
\begin{equation}\label{a5}
RP: H^0(\pp^2,\Omega^2_{\pp^2}(\mathcal{C}+m\mathcal{L})\otimes\mathcal{A}) \stackrel{\simeq}\longrightarrow H^0(\mathcal{C},\omega_{\mathcal{C}}(m\mathcal{L})\otimes \mathcal{A}_{\mathcal{C}}).
\end{equation}
\vskip2mm
\noindent
By \citep[Proposition p.33]{szpiro}, one has an isomorphism\footnote{One checks that the irreducibility assumption made in \citep[Proposition p.33]{szpiro} can be removed since the proof is local.}
\vskip2mm
\noindent
\[\begin{aligned}
Adj(d-3+m) & \stackrel{\simeq}\longrightarrow H^0(\pp^2,\Omega_{\pp^2}^2(\mathcal{C}+m\mathcal{L})\otimes \mathcal{A})\\
H & \longmapsto  \frac{Hdx\land dy}{F x^m}
\end{aligned}\]
\vskip2mm
\noindent
which combined with (\ref{a1}), (\ref{a4}) and (\ref{a5}) gives the isomorphism of Proposition \ref{p1}.
\end{proof}

\vskip2mm
\noindent
\begin{rem}
The isomorphism of Proposition \ref{p1} for $m=0$ and $\mathcal{C}$ irreducible is known as the Gorenstein Theorem \citep{gor} which asserts that the adjoint curves of degree $d-3$ cut out on $\mathcal{C}$ the complete canonical system of its normalization. For a nice down-to-earth presentation of adjoints and conductors, see \citep{ful}.
\end{rem}
\vskip2mm
\noindent
Recall from the introduction that one defines $A\subset k[y]$ to be the image of the projection
\[\begin{aligned}
Adj(d-2)&\longrightarrow  k[y]\\
H &\longmapsto  H(0,y).
\end{aligned}\]

\vskip3mm
\noindent
\begin{cor}\label{c0}
One has equality $dim_{k} \, A =d-\bar{s}$. 
\end{cor}
\vskip3mm
\noindent
\begin{proof}
If $H\in Adj(d-2)$ satisfies $H(0,y)\equiv 0$, then $H(x,y)=xH'(x,y)$ for some polynomial $H'$. Since the line $x=0$ does not contain any singularities of $\mathcal{C}$,  $H'$ is necessarily an adjoint of $F$ of degree $d-3$. In other words, one has an exact sequence of $k$-vector spaces
\vskip2mm
\noindent
\begin{equation}\label{exactadjoints}
0\longrightarrow Adj(d-3)  \longrightarrow  Adj(d-2) \longrightarrow A \longrightarrow 0
\end{equation}
\vskip2mm
\noindent
where the first map is the injective "multiplication by $x$" map and the second map is the restriction to $x=0$. 
It follows that 
\[\begin{aligned}
\dim(A) &=\dim Adj(d-2) - \dim Adj(d-3)\\
&= h^0(\omega_C(Z))-h^0(\omega_C)\\
&= d-\bar{s},
\end{aligned}\]
second equality using Proposition \ref{p1} and last equality using Proposition \ref{t2}.
\end{proof}

\vskip2mm
\noindent
The proof of Theorem \ref{adjmodx} follows.

\vskip3mm
\noindent

\begin{cor}\label{proofthm3} (Proof of Theorem $3$).
Given a basis of $Adj(d-2)$, one can compute a basis of $A$ within 
$$
\mathcal{O}((d-1)(g+d-\bar{s})(d-\bar{s})^{\omega-2})\subset \mathcal{O}(d^{\omega+1})
$$
arithmetic operations over $k$.
\end{cor}

\vskip3mm
\noindent
\begin{proof}
Consider the matrix $N$ whose set of rows is a basis of $Adj(d-2)$ evaluated at $x=0$, expressed in the natural basis of $k[y]$. So $N$ has $d-1$ columns and $g+d-\bar{s}$ rows (use Propositions \ref{t2} and \ref{p1}). By (\ref{exactadjoints}), one has $A=Im(N)$, and a basis of $A$ can be computed within the expected complexity 
\citep[Theorem 2.10]{stor}. The upper bound $\mathcal{O}(d^{\omega+1})$ follows from the well known inequality $g\le (d-1)(d-2)/2$.
\end{proof}

\vskip4mm
\noindent

\section{Recombinations follow}\label{S5} 

One has now all necessary information for solving recombinations. Let us consider first the rational case.

\vskip3mm
\noindent
\begin{cor}\label{c1}
One has an exact sequence of $k$-vector spaces
\vskip2mm
\noindent
\begin{equation*}
 0 \longrightarrow \langle \nu^{(1)},\ldots,\nu^{(s)} \rangle \longrightarrow k^n \stackrel{T}\longrightarrow A^{\vee}
\end{equation*}
\vskip2mm
\noindent
where $T$ sends $\nu=(\nu_1,\ldots,\nu_n)$ to the linear map

$$
H\longmapsto \sum_{i=1}^n \nu_i \,\, \Bigg(\sum_{f_i(\phi)=0}\frac{H(\phi)}{\partial_y F(0,\phi)}\Bigg). 
$$
\end{cor}

\vskip2mm
\noindent

\begin{proof}
By Lemma \ref{lemW}, one has 
\begin{equation*}
\langle \nu^{(1)},\ldots,\nu^{(s)} \rangle = Im(\alpha)\cap k^n
\end{equation*}
where one identifies $k^n\subset H^0(\mathcal{O}_Z)$ with the subspace of locally constant $k$-valued functions on $Z$. Proposition  \ref{t2} induces equality 
\begin{equation*}
Im(\alpha)\cap k^n =\Big\{\nu\in k^n, \quad\sum_{i=1}^n \nu_i res_{p_i}(\psi)=0 \quad \forall \,\,\psi \in H^0(\omega_C(Z))\Big\}.
\end{equation*}
\vskip1mm
\noindent
Let us compute the involved residues. By Proposition \ref{p1}, $\psi \in H^0(\omega_C(Z))$ is equal to
\begin{equation*}
\psi=\pi^*\Big(\frac{Hdx}{\partial_y F x}\Big)_{|C}
\end{equation*}
for a unique $H\in Adj(d-2)$. Let $\widehat{\mathcal{O}}_{C,p_i}$ be the completion of the regular local ring $\mathcal{O}_{C,p_i}$ with respect to its maximal ideal associated to $p_i$. The residue field  of $C$ at $p_i$ is equal to
$$
k_{p_i}=\frac{k[y]}{(f_i)}.
$$
The map $\pi$ being an isomorphism in a neighborhood of $p_i$, one has an isomorphism
\[\begin{aligned}
\widehat{\mathcal{O}}_{C,p_i} & \stackrel{\simeq}\longrightarrow  k_{p_i}[[t]]\\
\pi^* x            & \longmapsto   t\\
\pi^* y            & \longmapsto   a(t)
\end{aligned}\]
where $a\in k_{p_i}[[t]]$ is the unique series such that $a(0)$ is the residue class of $y$ in  $k_{p_i}$ and $F(t,a(t))\equiv 0$. In such a local system of coordinates, $\psi$ is equal to
$$
\psi = \frac{H(t,a(t))}{\partial_y F(t,a(t))}\frac{dt}{t}
$$
and it follows from the definition of residues that

\begin{equation*}\label{c13}
res_{p_i}(\psi)=Tr_{p_i} \Big(\frac{H(0,a(0))}{\partial_y F(0,a(0))}\Big) =\sum_{f_i(\phi)=0}\Big(\frac{H(0,\phi)}{\partial_y F(0,\phi)}\Big).
\end{equation*}
\vskip2mm
\noindent
Corollary \ref{c1} follows.
\end{proof}

\vskip1mm
\noindent
\begin{rem}
One always has $(1,\ldots,1)\in ker(T)$. This is nothing else than the Lagrange interpolation formula.
\end{rem}
\vskip1mm
\noindent
Let us now consider the absolute case.

\vskip1mm
\noindent

\begin{cor}\label{c2}
One has an exact sequence of $\bar{k}$-vector spaces
\begin{equation*}
 0 \longrightarrow \langle \bar{\nu}^{(1)},\ldots,\bar{\nu}^{(\bar{s})}\rangle\longrightarrow \bar{k}^d \stackrel{\bar{T}}\longrightarrow A^{\vee}\otimes_{k} \bar{k} \longrightarrow 0
\end{equation*}
\vskip2mm
\noindent
where $\bar{T}$ sends $\bar{\nu}=(\bar{\nu}_1,\ldots,\bar{\nu}_d)$ to the linear form
$$
H \longmapsto \sum_{i=1}^d \bar{\nu}_i \frac{H(\phi_i)}{\partial_y F(0,\phi_i)}.
$$
\end{cor}

\vskip3mm
\noindent

\begin{proof}
Apply Proposition  \ref{t2} and  repeat the proof of Corollary \ref{c1} over $\bar{k}$, with the curve $\bar{C}$ replacing $C$. Surjectivity of $\bar{T}$  follows from Corollary \ref{c0}. 
\end{proof}

\vskip3mm
\noindent

\begin{rem}\label{rk}In \citep{CL, L}, the authors solve recombinations using a system of $\mathcal{O}(d^2)$ equations. Corollary \ref{c1} and Corollary \ref{c2} give a much smaller number $d-\bar{s}$ of equations for recombinations. Moreover, the map $\bar{T}$ being surjective, $d-\bar{s}$ is the expected minimal number of linear conditions for recombinations in the absolute case. 
\end{rem}

\vskip2mm
\noindent
\section{Proofs of Theorems \ref{t1}, \ref{absolute} and \ref{absoluteproba}.}\label{S6}

In all what follows, one assumes that fast Fourier transform is used for polynomial multiplication, so that two univariate polynomials over $k$ of degree $\le m$ can be multiplied in softly linear time $\wt{\mathcal{O}}(m)$.

\vskip2mm
\noindent

\subsection{Proof of Theorem \ref{t1}.} One obtains the following algorithm.

\vskip4mm
\noindent
\textbf{Algorithm $1$ (deterministic rational factorization)}
\vskip2mm
\noindent

\textbf{Input :} $F\in k[x,y]$ that satisfies hypothesis (H).

\textbf{Output :} The rational factorization of $F$.
\vskip2mm
\noindent
\begin{itemize}
\item Step $1$. Compute a basis of $A$.
\item Step $2$. If $\dim A=1$, $F$ is irreducible. Otherwise, compute the irreducible factors $f_1,\ldots,f_n$ of $F(0,y)$ over $k$.
\item Step $3$. If $n=1$, $F$ is irreducible. Otherwise, build the matrix $M$ of the map $T$ of Corollary \ref{c1} by using Newton identities.
\item Step $4$. Compute the reduced echelon normal basis of $ker(M)$. One obtains the recombination vectors $\nu^{(1)},\ldots,\nu^{(s)}$.
\item Step $5$. Compute the factorization of $F(0,y)$ induced by the recombination vectors and lift it to the rational factorization of $F$.
\end{itemize}

\vskip4mm
\noindent

\begin{prop}\label{algo1}(Proof of Theorem \ref{t1}.)
Algorithm $1$ is deterministic and correct. Steps $3$, $4$ and $5$ take at most  
$$
\mathcal{O}(n(d-\bar{s})(d-s)^{\omega-2} + d^2)\subset \mathcal{O}(d^{\omega})
$$ 
arithmetic operations over $k$.
\end{prop}

\vskip4mm
\noindent

\begin{proof}
The algorithm is deterministic and correct thanks to Corollary \ref{c1}. Let us describe in more details the content and the complexity of steps $3$ to $5$.

\textit{Step $3$.} In order to build the matrix $M$, one has to compute 
$$
Tr_{k_{p_i}} \Big(\frac{H(y)}{\partial_y F(0,y)}\Big)
$$
for all $i=1,\ldots,n$ and for all $H$ running a basis of $A$. Inversion of $\partial_y F(0,y)$ and multiplication by $H$ in $k[y]/(f_i)$ take $\mathcal{O}(n_i)$ operations in $k$. Then $H/\partial_y F(0,y)\in k[y]/(f_i)$ is uniquely represented as a polynomial $a(y)=a_0+\cdots + a_{n_i-1}y^{n_i-1}$ with coefficients in $k$ and 
\begin{equation}\label{trace}
Tr_{k_{p_i}}\Big(\frac{H(y)}{\partial_y F(0,y)}\Big)=\sum_{j=0}^{n_i-1} a_j Tr_{k_{p_i}}(y^j).
\end{equation}
Thanks to the Newton identities, one can compute recursively the trace of $y^j$ from the traces of smaller powers of $y$ and from the coefficients of $f_i$ with $j$ multiplications and $j$ additions. So one computes traces of all involved powers of $y$ within  $\mathcal{O}(n_i^2)$ operations over $k$. Given these traces, and using (\ref{trace}), one computes the trace of $H/\partial_y F$ with $2n_i$ operations for each $H\in A$. By Corollary \ref{c0}, it follows that step $3$ costs $\sum_{i=1}^n\mathcal{O}(n_i^2+2n_i(d-\bar{s}))\subset \mathcal{O}(d^2)$ operations over $k$.

\textit{Step $4$.} The matrix $M$ has size $(d-\bar{s})\times n$ and rank $d-s$. One can compute the reduced echelon normal basis of the kernel of $M$ within $\mathcal{O}(n(d-\bar{s})(d-s)^{\omega-2})$ operations (\citep{stor}, Theorem 2.10). 

\textit{Step $5$.} Given a vector $\nu^{(j)}=(\nu^{(j)}_{i})\in \{0,1\}^n$ of the reduced echelon basis, one computes $F_j(0,y)=\prod f_i(y)^{\nu^{(j)}_{i}}$ for each rational irreducible factor $F_j$ of $F$. This requires $\wt{\mathcal{O}}(deg(F_j(0,y)))$ operations by the sub-product tree technique \citep[proof of Prop. 6]{L}, so a total cost of $\wt{\mathcal{O}}(d)$ operations. To compute the $F_j$'s, it's now enough to lift the induced equality $F(0,y)=F_1(0,y)\cdots F_s(0,y)$ modulo $(x)$ up to precision modulo $(x^{d+1})$. This costs $\wt{\mathcal{O}}(d^2)$ operations by using Newton quadratic iteration \citep[Theorem 15.18]{GG}. 
\end{proof}

\vskip3mm
\noindent

\subsection{Proofs of Theorems \ref{absolute} and \ref{absoluteproba}} In the absolute case, the delicate point is that Corollary \ref{c2} does not permit to solve recombinations with linear algebra over $k$. Moreover, it neither permits to describe the smallest finite extensions over which the irreducible absolute factors of $F$ are defined. To solve this problem, one rather relies our approach with the algorithms $8$ and $9$ in \citep{CL}, where the authors use the absolute partial fraction decomposition algorithm of Lazard-Rioboo-Trager \citep{LR}.
\vskip1mm
\noindent

Let $\phi$ be the residue class of $y$ in the ring $\mathbb{A}:=k[y]/(F(0,y))$. Any element $b\in \mathbb{A}$ can be uniquely represented as a finite sum
$$
b=\sum_{i=0}^{d-1} b_i \phi^{i}
$$
where $\coef(b,\phi^{i}):=b_i$ belongs to $k$. One introduces 

$$
L:=\Big\{v \in k^d,\quad \sum_{i=1}^d v_i \coef \Big(\frac{H(\phi)}{\partial_y F(0,\phi)}, \phi^{i-1}\Big) =0 \quad \forall \,\, H\in A\Big\}.
$$
\vskip2mm
\noindent
The vector space $L$ is related to the absolute recombinations by the following lemma. 

\vskip2mm
\noindent
\begin{lem}\label{imalpha} Let $V$ be the Vandermonde matrix of the roots $\phi_1,\ldots,\phi_d$ of $F(0,y)$. One has an isomorphism
$$
V^t : \langle \bar{\nu}^{(1)},\ldots,\bar{\nu}^{(\bar{s})}\rangle\stackrel{\simeq}{\longrightarrow} L\otimes_{k} \bar{k}.
$$
In particular, one has an isomorphism of $k$-vector spaces 
$$
B: Im(\alpha)\stackrel{\simeq}{\longrightarrow} L
$$
where $B=(Tr \phi^{i+j})_{i,j =0,\ldots,d-1}$, with $Tr:\mathbb{A}\to k$ the usual trace map.
\end{lem}

\vskip2mm
\noindent
\begin{proof}
One follows the proof of Proposition $4$ in \citep{CL}. Let $(v_1,\ldots,v_d)=V^t(w_1,\ldots,w_d)$ and let $b\in \mathbb{A}$. One has
\[\begin{aligned}
\sum_{i=1}^d v_i \coef (b,\phi^{i-1})
&=
\sum_{i=1}^d \big(\sum_{j=1}^d w_j \phi_j^i \big) \coef (b,\phi^{i-1})\\
&= 
\sum_{j=1}^d w_j \big(\sum_{i=1}^d \coef (b,\phi^{i-1})\phi_j^i\big)=\sum_{j=1}^d w_j b(\phi_j).
\end{aligned}\]
The first point then follows from Corollary \ref{c2} by taking $b=H(\phi)/\partial_y F(0,\phi)$. The second point follows from Lemma \ref{Wbar} since $V$ is the matrix of the evaluation map and $B=V^tV$ is the matrix of traces.
\end{proof}
\vskip2mm
\noindent
One can now rely on the factorization algorithms developed by Ch\`eze-Lecerf in the absolute case. One refers to their article \citep{CL} for details on the relations between absolute recombinations, absolute partial fraction decomposition, absolute Hensel lifting and absolute factorization.

\vskip6mm
\noindent
\textbf{Algorithm $2$ (deterministic absolute factorization).}
\vskip3mm
\noindent

\textbf{Input :} $F\in k[x,y]$ that satisfies hypothesis (H), with $k$ a field of characteristic $0$ or greater than $d(d-1)$.

\textbf{Output :} The absolute factorization of $F$.
\vskip2mm
\noindent
\begin{itemize}
\item Step $1$. Compute a basis of $A$.
\item Step $2$. Compute a basis of $L$.
\item Step $3$. Call Algorithm 8 in \citep{CL} with input $F$ and the basis of $L$.
\end{itemize}

\vskip4mm
\noindent

\begin{prop}\label{algo2}(Proof of Theorem \ref{absolute}.)
Algorithm $2$ is deterministic and correct. Steps $2$ and $3$ take at most  
$$
\wt{\mathcal{O}}(d(d-\bar{s})^{\omega-1}+\bar{s}d^3)\subset \wt{\mathcal{O}}(d^{4})
$$ 
arithmetic operations over $k$.
\end{prop}

\vskip4mm
\noindent

\begin{proof}
The algorithm is correct thanks to Lemma \ref{imalpha} combined with Proposition 4 p.15 and Theorem 5 p.15 in \citep{CL}. By definition, one has $L=ker(N)$, where the matrix $N$ is built from a basis of $A$ using one inversion in $A$ and $(d-\bar{s})$ multiplications in $\mathbb{A}$, so $\mathcal{O}(d(d-\bar{s}))$ operations over $k$. Then, computing a basis of $L=ker(N)$ requires $\mathcal{O}(d(d-\bar{s})^{\omega-1})$ operations over $k$. Finally, step $3$ costs $\wt{\mathcal{O}}(\bar{s}d^{3})$ operations over $k$ thanks to Proposition 10 p.24 in \citep{CL}.
\end{proof}
\vskip4mm
\noindent
The cost of Algorithm $2$ is dominated by the separation of residues in Algorithm 8 of \citep{CL} that ensures that the call to the Lazard-Rioboo-Trager algorithm returns a correct answer. If one rather deals with a random linear combination of the vectors of a basis of $L$, one obtains a probabilistic algorithm with smaller complexity.

\vskip6mm
\noindent
\textbf{Algorithm $3$ (probabilistic absolute factorization).}
\vskip3mm
\noindent

\textbf{Input :} $F\in k[x,y]$ that satisfies hypothesis (H), with $k$ a field of characteristic $0$ or greater than $d(d-1)$.

\textbf{Output :} The absolute factorization of $F$.
\vskip2mm
\noindent
\begin{itemize}
\item Step $1$. Compute a basis of $A$. 
\item Step $2$. Compute a basis of $L$.
\item Step $3$. Choose $c\in k^{\bar{s}}$, where $\bar{s}=d-dim (A)$.
\item Step $4$. Call Algorithm 9 in \citep{CL} with input $F$, the basis of $L$ and $c$.
\end{itemize}

\vskip4mm
\noindent

\begin{prop}\label{algo3}(Proof of Theorem \ref{absoluteproba}.)
Algorithm 3 either stops prematurely or return a correct answer. Moreover, there exists a polynomial 
$S\in\bar{k} [C_1,\ldots,C_{\bar{s}}]$ 
of degree at most $\bar{s}(\bar{s}-1)$ such that the answer is correct whenever $S(c)\ne 0$. In any cases, steps $2$, $3$ and $4$ take at most  
$$
\wt{\mathcal{O}}(d(d-\bar{s})^{\omega-1}+d^{\frac{\omega+3}{2}})
$$ 
arithmetic operations over $k$. 
\end{prop}

\vskip4mm
\noindent

\begin{proof}
The proof is the same as for Proposition \ref{algo2}, using now  Proposition 11 p.25 in \citep{CL}.
\end{proof}

\vskip4mm
\noindent

\begin{rem}
The restriction hypothesis on the characteristic of $k$ ensures the possibility to separate the residues and to apply a fast absolute multi-factor Hensel lifting in \citep[Algorithm 9]{CL}. 
\end{rem}
\vskip2mm
\noindent

\section{The case $F(0,y)$ non separable}\label{S7}

When $F(0,y)$ is not separable modulo $(x)$, one is tempted to choose another fiber $x=a$ for which $F(a,y)$ that satisfies hypothesis $(H)$. There are two main reasons to develop a recombination algorithm along a critical fiber. First, when the field $k$ has small positive characteristic, a regular fiber may not exist. Second, one shows here that working along a singular fiber may in fact be an opportunity to speed-up the algorithm. 

\vskip2mm
\noindent

In order to simplify, one supposes here that $k=\bar{k}$. One supposes too that $F(0,y)$ has degree $d$ (the general case follows easily by computing residues at $y=\infty$). 

\vskip2mm
\noindent

Our results generalise well to the non separable case, the main difference being related to the computation of residues. Let $Z=C\cap L$ where $L=\pi^*(\mathcal{L})$. In contrast to the previous sections, $L$ and $Z$ need not to be reduced anymore. The support of $Z$ consists now in $r\le d$ closed points $p_1,\ldots, p_r$ in one-to-one correspondance with the irreducible analytic branches of $\mathcal{C}$ along the line $\mathcal{L}$. The recombination vectors may now be defined in the smaller ambient space $k^r$, where $\mu^{(j)}\in k^r$ is defined to have $i$th coordinate equal to $1$ if $p_i\in C_j$ and equal to $0$ otherwise. By identifying $k^r$ with the vector subspace of $H^0(\mathcal{O}_Z)$ of locally constant functions on $Z$ with values in $k$ \textit{and with zero nilpotent part}, one obtains the analoguous of Lemma \ref{lemW} 
\begin{equation}\label{rec}
\langle \nu^{(1)},\ldots,\nu^{(s)} \rangle = Im(\alpha)\cap k^r,
\end{equation}
\vskip2mm
\noindent
where $\alpha:H^0(\mathcal{O}_C) \to H^0(\mathcal{O}_Z)$ still stands for the restriction map. One obtains the following generalization of Corollary \ref{c1} :
\vskip4mm
\noindent
\begin{prop}\label{pfin}
One has an exact sequence of $k$-vector spaces
\vskip2mm
\noindent
\begin{equation*}
 0 \longrightarrow \langle \nu^{(1)},\ldots,\nu^{(s)} \rangle \longrightarrow k^r \stackrel{T}\longrightarrow A^{\vee}
\end{equation*}
\vskip2mm
\noindent
where $T$ sends $\nu=(\nu_1,\ldots,\nu_r)$ to the linear map
$$
H\longmapsto \sum_{i=1}^{r}\nu_i \,res_{p_i}\Big(\pi^*\Big(\frac{H(y)dy}{F(0,y)}\Big)\Big). 
$$
\end{prop}
\vskip4mm
\noindent

\begin{proof}
The exact sequence of Proposition \ref{t2} remains valid in this new context. Combined with (\ref{rec}), one obtains that $\langle \nu^{(1)},\ldots,\nu^{(s)}\rangle = ker \,\wt{T}$ where 
$
\wt{T} : k^r \rightarrow H^0(\omega_{C}(Z))^{\vee}
$
sends $\nu=(\nu_1,\ldots,\nu_r)$ to the linear map
$$
\psi\longmapsto  \sum_{i=1}^{r}\nu_i res_{p_i}(\psi). 
$$
The main difference concerns the computation of residues. Let $\psi\in H^0(\omega_{C}(Z))$. By the proof of Proposition \ref{p1}, one has 
$$
\psi=RP_C(\Psi),\qquad \Psi = \pi^*\Bigg(\frac{H(x,y)dx\land dy}{xF(x,y)}\Bigg)\in H^0(\Omega^2_X(C+L))
$$
for a unique $H\in Adj(d-2)$, and where $RP_C$ stands for the Poincar\'e residue along $C$. Let $p\in C\cap L$. By \citep[Theorem 6.3 and Remark 6.9]{couvreur}, one obtains equality
$$
res_p(\psi)=res_p\Big(RP_C(\Psi))\Big)=res_p\Big(RP_L(\Psi))\Big).
$$
where the last residue stands for residue of $1$-form on $L$. Since $L=\pi^*(\mathcal{L})$ and the Poincar\'e residue commutes with the pull-back, one obtains equality
$$
RP_L\Bigg(\pi^*\Bigg(\frac{H(x,y)dx\land dy}{xF(x,y)}\Bigg)\Bigg)=\pi^*\Bigg(RP_{\mathcal{L}}\Bigg(\frac{H(x,y)dx\land dy}{xF(x,y)}\Bigg)\Bigg)=\pi^*\Bigg(\frac{H(0,y)dy}{F(0,y)}\Bigg).
$$
Proposition \ref{pfin} follows.
\end{proof}

So as soon as the curve $\mathcal{C}$ has a small number of irreducible branches intersecting $\mathcal{L}$, Proposition \ref{pfin} permits to solve the recombination problem in a smaller ambient space. The price to pay is that one can {\it a priori} not compute residues directly in $\pp^2$ as in Section \ref{S2}, but one may really need to compute residues in $X$, using local coordinates or Puiseux series. 
Nevertheless, in the important case of $\mathcal{C}$ locally irreducible at $(0,y_p)\in \mathcal{C}\cap \mathcal{L}$, there is exactly one point $p \in C$ such that $\pi(p)=(0,y_p)$ and the residue can be computed directly in $\pp^2$ :
$$
res_{p}\Bigg(\pi^*\Big(\frac{H(y)dy}{F(0,y)}\Big)\Bigg)=res_{y_p}\Big(\frac{H(y)dy}{F(0,y)}\Big).
$$
Of course $y_p$ may be now be a multiple root of $F(0,y)$ so the residue computation may involve higher order derivatives of $H$ and $F$. 

\vskip4mm
\noindent
\textbf{Example.} {\it Let us illustrate Proposition \ref{pfin} on a simple example. Suppose that one wants to factorize
$$
F(x,y)=y^5+y^4-xy^3 - y^3 -2xy^2 -y^2 + x^2 + xy +x.
$$
over $\qp$. The irreducible factorization of $F$ mod $(x)$ is
$$
F(0,y)=y^2(y+1)^2(y-1),
$$ 
with two double roots $-1$ and $0$ and one simple root $1$. Since $\partial_x F(0,-1)$ and $\partial_x F(0,0)$ do not vanish, the curve $\mathcal{C}$ is smooth and tangent to $\mathcal{L}$ at these two points, and transversal to $\mathcal{L}$ at $(0,1)$. In particular, $\mathcal{C}$ has only $3$ irreducible branches intersecting $\mathcal{L}$ at distinct points and the recombinations will hold in the ambient space $k^3$ rather than in the bigger space $k^5$ inherent to a choice of a regular fiber.
\vskip2mm
\noindent
One has here
$$
Adj(d-2)=\langle y^3-y-1,y^2-x,y^2-x,xy^2-x^2,y^3-xy\rangle
$$
from which it follows that 
$$
A=\langle y+1,y^2,y^3 \rangle.
$$
\vskip2mm
\noindent
By Lemma \ref{lemDim}, the curve $\mathcal{C}$ has $deg(F)-dim(A)=2$ absolute irreducible components. Let $H \in A$. Since $F$ is locally irreducible at $(0,0)$, one has equality  

\[\begin{aligned}
res_{0}\Big(\frac{H(y)dy}{F(0,y)}\Big)=res_{0}\Big(\frac{H(y)dy}{y^2(y+1)^2(y-1)}\Big)= H'(0)+H(0).
\end{aligned}\]
\vskip2mm
\noindent
In the same way, a simple calculation gives  

\[\begin{aligned}
res_{-1}\Big(\frac{H(y)dy}{F(0,y)}\Big)=\frac{-2H'(-1)-5H(-1)}{4}\qquad {\rm and}\qquad 
res_{1}\Big(\frac{H(y)dy}{F(0,y)}\Big)&=&\frac{H(1)}{4}.
\end{aligned}\]
\vskip2mm
\noindent
One deduces that the $3\times 3$ matrix of the map $T$ in Proposition \ref{pfin} is
$$
M=\begin{pmatrix}
-1/2 & -1/4 & -1/4  \\ 
0 & 0 & 0 \\ 
1/2 & 1/4 & 1/4 \\
\end{pmatrix}
$$
so that $ker(M)=((0,1,0),(1,0,1))$. One deduces the irreducible rational decomposition 
$$
\mathcal{C}=\mathcal{C}_1\cup \mathcal{C}_2
$$
\vskip2mm  
\noindent
where $\mathcal{C}_1$ is a conic tangent to $\mathcal{L}$ at $(0,0)$ and $\mathcal{C}_2$ is a (smooth) 
cubic tangent to $\mathcal{L}$ at $(0,1)$ and transversal to $\mathcal{L}$ at $(0,-1)$. 

\vskip3mm  
\noindent
\begin{equation*}
\begin{array}{c}
\resizebox{5.cm}{!}{\includegraphics{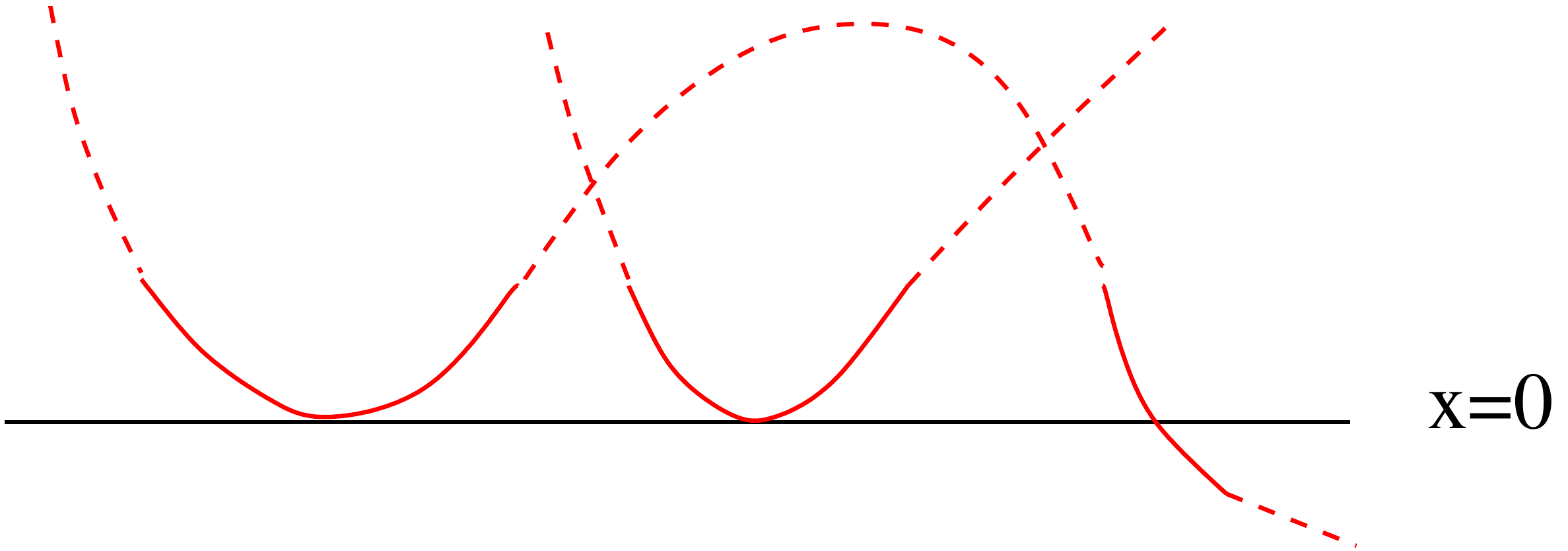}} 
\end{array}
\end{equation*}
\vskip3mm  
\noindent
Since the induced factors of $F$ are coprime modulo $(x)$, one can compute them use a multifactor Hensel lifting \citep[Algorithm 15.17]{GG} up to a sufficiently high precision (mod $(x^3)$ in our case). One obtains finally
$$
F(x,y)=(y^2-x)((y+1)^2(y-1)-x).
$$}
\vskip3mm
\noindent

In summarize, one has shown here that working along a critical fiber may be an opportunity to speed-up the algorithm, at least when $F$ satisfies the weaker hypothesis
\begin{equation*}
{\rm (H')} \qquad \mathcal{C} \,\,{\rm  is \,\,analytically \,\,irreducible\,\, at\,\, each \,\,point \,\,of\,\,} \mathcal{C}\cap \mathcal{L}.
\end{equation*}
First, the univariate factorization of $F(0,y)$ is faster since it is reduced to a fast separable factorization \citep{L2} plus some univariate factorizations of smaller degrees. Second, recombinations are faster since they hold in a smaller ambient space of dimension the number of distinct roots of $F(0,y)$ (or irreducible factors in the rational case). This fact is well illustrated in a previous work of the author \citep{W2} who developed a lifting and recombination algorithm based on the toric resolution of the singularities of $\mathcal{C}$ along the line at infinity. 

\vskip3mm
\noindent

\begin{rem}
One can show that under hypothesis $(H')$, building the matrix of the map $T$ has the same cost in the separable and non separable cases. The fact that the computations of the residues may involve $f_i$-adic expansions with higher precision for each irreducible factors $f_i$ of $F(0,y)$ is compensed by the fact that the sum of the degrees of the $f_i$'s decreases. 
\end{rem}
\vskip3mm
\noindent

\begin{rem}
Although  Proposition \ref{pfin} still permit to solve recombinations even when $F$ does not satisfy $(H')$, the problem resides in the fact that the irreducible factors of $F$ may not be coprime modulo $(x)$ and can not be computed with Hensel's lemma. This problem will be explored in a further work. 
\end{rem}

\vskip2mm
\noindent

\section{Don't touch the cusps}\label{S8}

It turns out that the computation of $A$ does not necessarily require to compute the all resolution of singularities of $\mathcal{C}$. 
Namely, let us consider the factorization of $\pi$
$$
X\longrightarrow X_0 \stackrel{\pi_0}\longrightarrow \pp^2
$$ 
where $\pi_0$ is the minimal composition of blow-ups under which the strict transform $C_0$ of $\mathcal{C}$ is every where locally irreducible. 
Then one can check that all our results (Lemmas \ref{lemW}, \ref{lemDim}, \ref{Wbar} and the key Proposition \ref{p1}) remain valid with $C_0$ replacing $C$ and with the arithmetic genus $p_a(C_0)\ge g$ of $C_0$ replacing the geometric genus of $\mathcal{C}$ (the proofs mainly only use that the irreducible and connected components of $C$ coincide). Then, one checks easily that one has an exact sequence
$$
0 \longrightarrow Adj_0(d-3)\stackrel{\times x}\longrightarrow Adj_0(d-2) \stackrel{x=0}\longrightarrow A \longrightarrow 0
$$
\vskip1mm
\noindent
where $Adj_0(k)$ is defined similarly as $Adj(k)$ with the map $\pi_0$ replacing $\pi$ in Definition \ref{defadj}. For instance, one needs not to desingularize cusps for computing $A$, which is of course natural from our factorization point of view. Note that there are easy local irreducibility sufficient criterions that can be directly read off from the Newton polygon of the singularity (for cusps for instance).

\vskip1mm
\noindent

\section{Conclusion}\label{S9}

One has established the bridge between locally constant functions \citep{duval, ragot}
and lifting and recombinations algorithms \citep{L, CL}. One has shown that the computation of adjoint polynomials allows to solve the recombinations problems without using Hensel lifting and with the expected number of linear equations. Although one believes that our approach uses too strongly the geometry of singularities, it may be useful for some particular polynomials, for instance if the irreducible components are smooth and intersect in few points. Moreover, one has discussed the possibility of speed-up the algorithm when $F$ is not separable modulo $(x)$. 
In a further work, one plans to develop some intermediary algorithms based only on the resolution of some of the singularities, in the vein of \citep{W2} that uses the toric resolution of the singularities at infinity. In any cases, the power of using singularities for factorization depends strongly on complexity issues in the algorithmic theory of singularities, especially on the Newton-Puiseux algorithm.

\vskip3mm
\noindent

\bibliographystyle{elsart-harv}


\end{document}